\documentclass[12pt,a4paper]{amsart}
\usepackage[T1]{fontenc}
\usepackage{mathpazo}
\usepackage[margin=30mm,footskip=15mm]{geometry}
\usepackage{microtype,amssymb,mathtools}
\usepackage{booktabs,tabularx,array,needspace}
\usepackage{xcolor,tikz}
\usepackage{float}
\usetikzlibrary{arrows.meta}
\definecolor{maroon}{RGB}{128,0,0}
\definecolor{pathblue}{RGB}{24,65,105}
\usepackage[colorlinks=true,linkcolor=pathblue,citecolor=pathblue,urlcolor=pathblue]{hyperref}
\hypersetup{pdftitle={Associahedra are 8-colourable},pdfauthor={David R. Wood},
 pdfsubject={An 8-colouring of associahedron graphs},
 pdfkeywords={associahedron, graph colouring, tree rotation}}
\newtheorem{theorem}{Theorem}
\newtheorem{lemma}[theorem]{Lemma}
\newcommand{\F}{\mathbb F_2}
\newcommand{\An}{\mathcal A_n}
\newcommand{\lab}{\lambda}
\newcommand{\defterm}[1]{\textcolor{maroon}{\emph{#1}}}
\newcommand{\defmath}[1]{{\color{maroon}#1}}
\title{Associahedra are 8-colourable}
\author{David R. Wood}
\newcommand{\paperdate}{12 November 2026}
\date{}
\makeatletter
\renewcommand{\@setauthors}{%
  \begingroup
  \centering
  \vskip12pt
  {\normalsize\scshape\authors\par}
  \vskip4pt
  {\normalsize School of Mathematics, Monash University\par
  Melbourne, Australia\par}
  \vskip4pt
  {\small\href{mailto:david.wood@monash.edu}{\texttt{david.wood@monash.edu}}\par}
  \vskip6pt
  {\normalsize\paperdate\par}
  \endgroup
}
\makeatother
\begin{document}
\begin{abstract}
Determining the chromatic number of associahedra is a long-standing open
problem. Tan [\href{https://arxiv.org/abs/2611.05323}{arXiv:2611.05323}]
recently proved that associahedra are 12-colourable, which is the first
constant upper bound. This paper adapts Tan's argument to show
that associahedra are 8-colourable.

\smallskip\noindent\textbf{AI Statement.} The construction and proof were developed in
collaboration with OpenAI Codex. The author checked the proof and thoroughly
revised the paper.
\end{abstract}
\maketitle
\vspace{-2\baselineskip}

\section{Introduction}

Associahedra arise in many branches of mathematics. In algebraic
topology, Stasheff introduced associahedra to describe associativity up
to homotopy~\cite{Stasheff1963}. In algebra and combinatorics, the
vertices and edges of an associahedron form the Hasse diagram of the
Tamari lattice of bracketings~\cite{Tamari1962}. In discrete geometry,
associahedra are polytopes with many different
realisations~\cite{Loday2004}. In the theory of cluster algebras,
associahedra are the generalised associahedra of type~$A$~\cite{FZ2003}.
In algebraic geometry, copies of an associahedron tile the real points
of the moduli space of stable genus-zero curves with marked
points~\cite{Devadoss1999}. In physics, associahedra appear in the
positive geometry of scattering amplitudes~\cite{ABHY2018}. In computer
science, the graph of an associahedron is the rotation graph of binary
trees, which models the restructuring of binary search
trees~\cite{STT1988}. In computational geometry, the same graph is the
flip graph of the triangulations of a convex polygon.

\begin{figure}[htbp]
\centering
\input{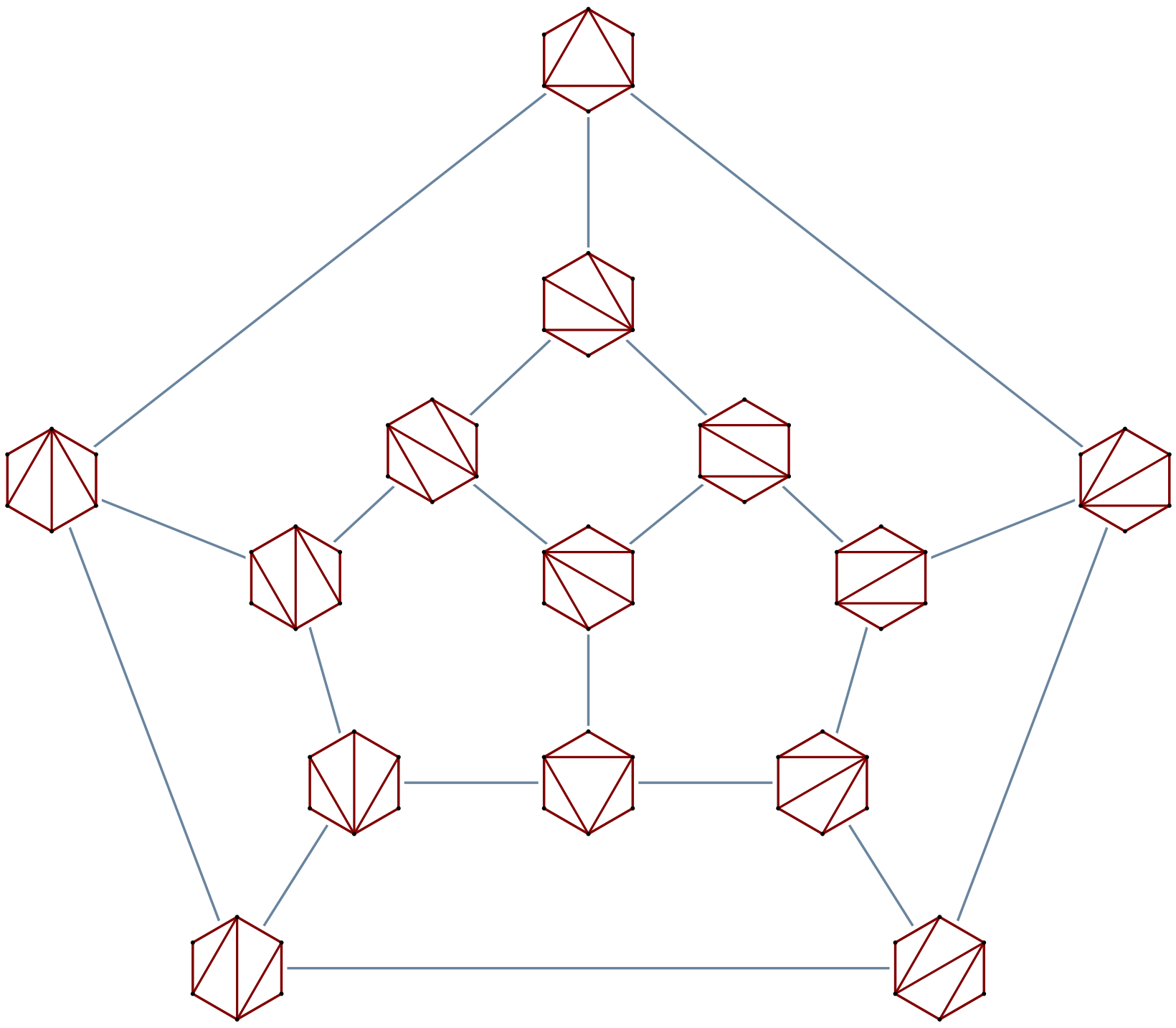}
\caption{The associahedron graph $\mathcal A_6$.}
\label{fig:associahedron6}
\end{figure}

A \defterm{triangulation} of a convex polygon is a maximal set of pairwise
noncrossing diagonals. For each integer $n\geqslant3$, the
\defterm{associahedron graph} $\defmath{\An}$ is the graph whose vertices
are the triangulations of a convex $n$-gon. Two triangulations are
adjacent in $\An$ if one of them can be obtained from the other by
replacing a single diagonal. Figure~\ref{fig:associahedron6} shows
$\mathcal A_6$. A \defterm{colouring} of a graph assigns a
colour to each vertex so that adjacent vertices receive different
colours. The \defterm{chromatic number} $\defmath{\chi(G)}$ of a graph
$G$ is the minimum number of colours in a colouring of $G$.

Fabila-Monroy, Flores-Pe\~naloza, Huemer, Hurtado, Urrutia and
Wood~\cite{Fabila2009} initiated the study of $\chi(\An)$.%
\begin{table}[!b]
\caption{Upper bounds on $\chi(\An)$. The date refers to the first release.}
\label{tab:history}
\centering
\small
\setlength{\tabcolsep}{4pt}
\renewcommand{\arraystretch}{1.13}
\begin{tabularx}{\textwidth}{@{}p{28mm}>{\raggedright\arraybackslash}X p{46mm}@{}}
\toprule
Date & Reference & Upper bound\\
\midrule
March 2007 & Fabila-Monroy et al.~\cite{Fabila2009}
 & $\lceil n/2\rceil$ and $O(n/\log n)$\\
November 2018 & Addario-Berry, Reed, Scott and Wood~\cite{ARSW2026}
 & $O(\log n)$\\
September 2026 & Berendsohn, Cardinal, Iacono and Kozma~\cite{BCIK2026}
 & $6\log_2(n-2)+18$\\
September 2026 & Oum and Wood~\cite{OW2026}
 & $10\log_2\log_2 n$ for $n\geqslant4$\\
November 2026 & Tan~\cite{Tan2026} & $12$\\
November 2026 & Theorem~\ref{thm:main} & $8$\\
\bottomrule
\end{tabularx}
\end{table}
These authors gave an explicit colouring of $\An$ with
$\lceil n/2\rceil$ colours. They also deduced the upper bound
$\chi(\An)=O(n/\log n)$ from Johansson's colouring theorem for
triangle-free graphs~\cite{Johansson1996}. In addition, Fabila-Monroy et
al.\ made two conjectures. The first conjecture states that
$\chi(\An)=O(\log n)$. The second conjecture states that $\chi(\An)$ is
unbounded. Addario-Berry, Reed, Scott and Wood~\cite{ARSW2026} proved the
first conjecture. Later, Berendsohn, Cardinal, Iacono and
Kozma~\cite{BCIK2026} gave a simpler colouring with $O(\log n)$
colours. Oum and Wood~\cite{OW2026} then improved the upper bound to
$O(\log\log n)$. Recently, Tan~\cite{Tan2026} disproved the second
conjecture of Fabila-Monroy et al.\ by showing that $\An$ has a
$12$-colouring for all $n\geqslant3$. Table~\ref{tab:history} lists these
upper bounds on $\chi(\An)$.

We establish the following improved bound.
\begin{theorem}\label{thm:main}
$\chi(\An)\leqslant8$ for every integer $n\geqslant3$.
\end{theorem}

The best known lower bound is $\chi(\mathcal A_{10})\geqslant4$, due to
Ruy Fabila-Monroy~\cite{ARSW2026}. It is open whether
$\chi(\An)\geqslant5$ for some $n$~\cite{Tan2026}. Oum and
Wood~\cite{OW2026} verified by computer that
$\chi(\mathcal A_{15})\leqslant4$. Since $\An$ is isomorphic to an induced
subgraph of $\mathcal A_{n+1}$, it follows that $\chi(\An)=4$ for each
$n\in\{10,\dots,15\}$.

Section~\ref{sec:proof} proves Theorem~\ref{thm:main}.
Section~\ref{sec:lower} considers colourings of the type that are used in
the proof of Theorem~\ref{thm:main}, and shows that the set of colours has
at least $8$ elements.

\section{The Proof}\label{sec:proof}

\begin{figure}[!h]
\centering
\begin{tikzpicture}[x=1cm,y=1cm,
 vertex/.style={circle,fill=black,inner sep=1.7pt},
 every node/.style={font=\small},line width=.65pt]
 \coordinate (r) at (0,0);
 \coordinate (s) at (-.85,-.85);
 \coordinate (u) at (-1.6,-1.75);
 \coordinate (v) at (-.1,-1.75);
 \coordinate (w) at (1.3,-1.75);
 \draw (r)--(s) (r)--(w) (s)--(u) (s)--(v);
 \node[vertex] at (r) {};
 \node[vertex] at (s) {};
 \foreach \root/\letter in {u/U,v/V,w/W}{
   \draw[fill=black!7] (\root)--++(-.45,-.7)--++(.9,0)--cycle;
   \node[vertex] at (\root) {};
   \node at ([yshift=-4.7mm]\root) {$\letter$};
 }
 \node at (0,-2.8) {$((U,V),W)$};
 \draw[{Stealth[length=2mm]}-{Stealth[length=2mm]}]
   (1.95,-.8)--(3.05,-.8);
 \coordinate (rp) at (5.2,0);
 \coordinate (t) at (6.05,-.85);
 \coordinate (up) at (3.9,-1.75);
 \coordinate (vp) at (5.3,-1.75);
 \coordinate (wp) at (6.8,-1.75);
 \draw (rp)--(up) (rp)--(t) (t)--(vp) (t)--(wp);
 \node[vertex] at (rp) {};
 \node[vertex] at (t) {};
 \foreach \root/\letter in {up/U,vp/V,wp/W}{
   \draw[fill=black!7] (\root)--++(-.45,-.7)--++(.9,0)--cycle;
   \node[vertex] at (\root) {};
   \node at ([yshift=-4.7mm]\root) {$\letter$};
 }
 \node at (5.2,-2.8) {$(U,(V,W))$};
\end{tikzpicture}
\caption{A rotation of an ordered full binary tree.}
\label{fig:rotation-schematic}
\end{figure}
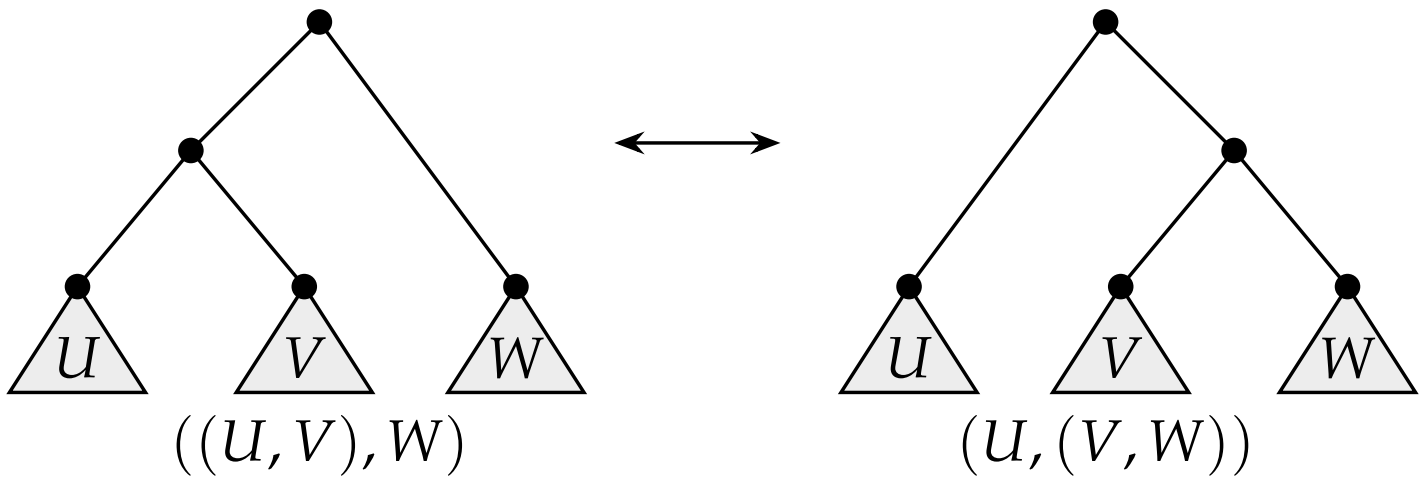

We begin by representing triangulations by trees.
An \defterm{ordered full binary tree} is a finite rooted tree in which
every vertex has either $0$ or $2$ children. A vertex with $0$ children
is a \defterm{leaf}. A vertex with $2$ children is \defterm{internal}.
The two children of each internal vertex are designated as its left
child and its right child. A \defterm{left edge} joins a parent to its
left child. A \defterm{right edge} joins a parent to its right child.
We record the direction of a left edge by $\defmath{L}$ and the direction
of a right edge by $\defmath{R}$. For a vertex $v$ of a tree $T$, the
\defterm{rooted subtree} of $T$ at $v$ consists of $v$ and all
descendants of $v$ in $T$. For trees $U$ and $V$, let $\defmath{(U,V)}$
denote the tree that is obtained by adding a new root whose left child is
the root of $U$ and whose right child is the root of $V$.
A \defterm{rotation} replaces a rooted subtree of the form $((U,V),W)$ by
$(U,(V,W))$, or replaces a rooted subtree of the form $(U,(V,W))$ by
$((U,V),W)$. Figure~\ref{fig:rotation-schematic} illustrates a rotation.

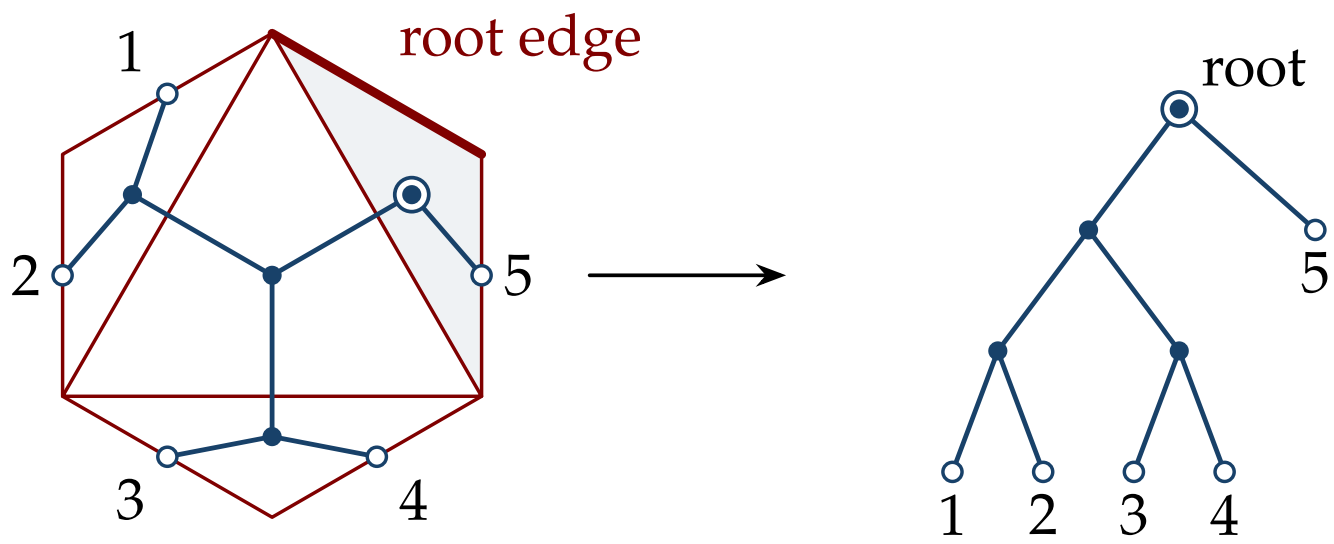
\begin{figure}[!h]
\centering
\begin{tikzpicture}[x=1cm,y=1cm,line cap=round,line join=round,
 every node/.style={font=\small}]
\coordinate (h0) at (0,1.6);
\coordinate (h1) at (-1.38564,.8);
\coordinate (h2) at (-1.38564,-.8);
\coordinate (h3) at (0,-1.6);
\coordinate (h4) at (1.38564,-.8);
\coordinate (h5) at (1.38564,.8);
\fill[pathblue!7] (h0)--(h4)--(h5)--cycle;
\draw[maroon,line width=.65pt] (h0)--(h1)--(h2)--(h3)--(h4)--(h5)--cycle;
\draw[maroon,line width=.65pt] (h0)--(h2)--(h4)--(h0);
\draw[maroon,line width=1.7pt] (h0)--(h5);
\node[above right,text=maroon] at (.69282,1.2) {root edge};
\coordinate (a) at (.92376,.53333);
\coordinate (b) at (0,0);
\coordinate (c) at (-.92376,.53333);
\coordinate (d) at (0,-1.06667);
\coordinate (l1) at (-.69282,1.2);
\coordinate (l2) at (-1.38564,0);
\coordinate (l3) at (-.69282,-1.2);
\coordinate (l4) at (.69282,-1.2);
\coordinate (l5) at (1.38564,0);
\draw[pathblue,line width=.9pt] (a)--(b) (a)--(l5)
 (b)--(c) (b)--(d) (c)--(l1) (c)--(l2) (d)--(l3) (d)--(l4);
\foreach \v in {a,b,c,d}{\fill[pathblue] (\v) circle[radius=1.8pt];}
\draw[draw=pathblue,fill=white,line width=.7pt] (a) circle[radius=3.2pt];
\fill[pathblue] (a) circle[radius=1.8pt];
\foreach \v in {1,...,5}{\draw[draw=pathblue,fill=white,line width=.65pt] (l\v) circle[radius=1.8pt];}
\node[above left] at (l1) {$1$};
\node[left] at (l2) {$2$};
\node[below left] at (l3) {$3$};
\node[below right] at (l4) {$4$};
\node[right] at (l5) {$5$};
\draw[-{Stealth[length=2mm]},line width=.7pt] (2.1,0)--(3.4,0);
\begin{scope}[shift={(5.7,0)}]
\coordinate (a) at (.3,1.1);
\coordinate (b) at (-.3,.3);
\coordinate (c) at (-.9,-.5);
\coordinate (d) at (.3,-.5);
\coordinate (l1) at (-1.2,-1.3);
\coordinate (l2) at (-.6,-1.3);
\coordinate (l3) at (0,-1.3);
\coordinate (l4) at (.6,-1.3);
\coordinate (l5) at (1.2,.3);
\draw[pathblue,line width=.9pt] (a)--(b) (a)--(l5)
 (b)--(c) (b)--(d) (c)--(l1) (c)--(l2) (d)--(l3) (d)--(l4);
\foreach \v in {a,b,c,d}{\fill[pathblue] (\v) circle[radius=1.8pt];}
\draw[draw=pathblue,fill=white,line width=.7pt] (a) circle[radius=3.2pt];
\fill[pathblue] (a) circle[radius=1.8pt];
\node[above right] at (a) {root};
\foreach \v in {1,...,5}{
 \draw[draw=pathblue,fill=white,line width=.65pt] (l\v) circle[radius=1.8pt];
 \node[below] at (l\v) {$\v$};
}
\end{scope}
\end{tikzpicture}
\caption{The tree--triangulation bijection.}
\label{fig:tree-triangulation-bijection}
\end{figure}

The classical tree--triangulation bijection identifies the vertices of
$\An$ with the ordered full binary trees that have $n-1$ leaves. Under
this bijection, two vertices of $\An$ are adjacent if and only if the
corresponding trees differ by one rotation. Sleator, Tarjan and
Thurston~\cite[Section~2.2, Lemma~1]{STT1988} give an explicit
construction of this bijection, as illustrated in
Figure~\ref{fig:tree-triangulation-bijection}.

\begin{figure}[!t]
\centering
\input{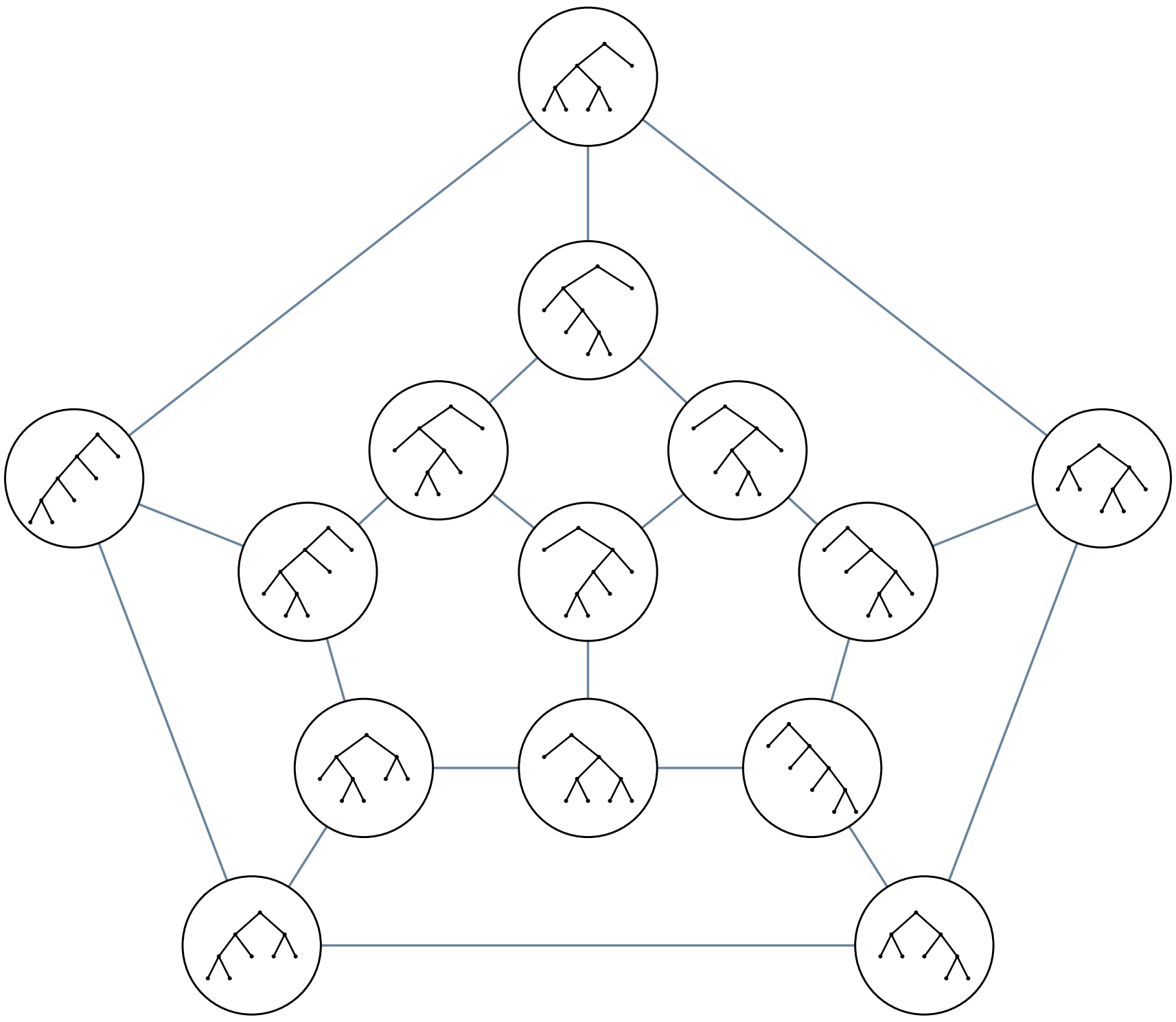}
\caption{The associahedron graph $\mathcal A_6$ represented by ordered
full binary trees with five leaves. Each tree occupies the same position
as its corresponding triangulation in Figure~\ref{fig:associahedron6}.}
\label{fig:associahedron6-trees}
\end{figure}

From now on, every tree is an ordered
full binary tree, and we work with trees instead of triangulations.
Figure~\ref{fig:associahedron6-trees} shows the trees corresponding to
the triangulations in Figure~\ref{fig:associahedron6}.

We now define a recursive labelling of the vertices of a tree, following
Tan~\cite{Tan2026}, who introduced such a labelling with labels in
$\mathbb Z_3$. Our recursion uses four labels instead of three.
Let $\defmath{\F}\coloneqq\{0,1\}$ be the field of integers modulo $2$.
Define a binary operation $\defmath{Q}$ on $\F^2$ by
\begin{equation*}
 \defmath{Q(a,b)}\coloneqq \bigl(a_0+a_1+b_0,\ a_1+b_0+b_1+1\bigr),
 \qquad a=(a_0,a_1),\quad b=(b_0,b_1).
\end{equation*}
Every leaf receives the label $(0,0)$. Consider an internal vertex whose
left child has label $a$ and whose right child has label $b$. This
internal vertex receives the label $Q(a,b)$. For a tree $T$, let
$\defmath{\lab(T)}$ denote the label of the root of $T$. By construction,
the label of a vertex depends only on the rooted subtree at that vertex.

We next use the labels to select one child at each internal vertex.
Consider an internal vertex whose left child has label $a$ and whose
right child has label $b$. Define
\begin{equation*}
 \defmath{D(a,b)}\coloneqq a_1+b_0.
\end{equation*}
If $D(a,b)=0$, then the edge from this vertex to its left child is
\defterm{preferred}. If $D(a,b)=1$, then the edge from this vertex to its
right child is preferred. In both cases, the other edge from this vertex
to a child is \defterm{unpreferred}. The child at the lower end of the
preferred edge is the \defterm{preferred child}. Thus every internal
vertex has exactly one preferred child.

Let $T$ be a tree. Consider the subgraph of $T$ that contains every
vertex of $T$ and exactly the preferred edges of $T$. In this subgraph,
each vertex is incident with at most one edge to a child and at most
one edge to its parent. Since $T$ has no cycles, every component of this subgraph is a
path. These components are the \defterm{preferred paths} of $T$. A
preferred path may consist of a single vertex. Each preferred path has a
unique vertex that is nearest to the root of $T$. This vertex is the
\defterm{top} of the preferred path. We orient each preferred path away
from its top. Since every internal vertex has a preferred child, every
preferred path descends from its top to a leaf. The terms preferred child
and preferred path are borrowed from link-cut trees~\cite{ST1983} and
Tango trees~\cite{Tango2007}. This method of forming
paths, in which one child is selected at each internal vertex and every
path ends at a leaf, follows Oum and Wood~\cite[Section~3]{OW2026}.
Tan~\cite{Tan2026} uses the recursive labels to determine a set of
paths, called highlighted paths. Our labels instead select exactly one
child at each internal vertex.

The \defterm{word} of a preferred path is the sequence of directions
$L$ and $R$ of its edges, listed in order from its top. Oum and
Wood~\cite[Section~3]{OW2026} also record paths by words in $L$ and $R$.

\begin{figure}[!htbp]
\centering
\input{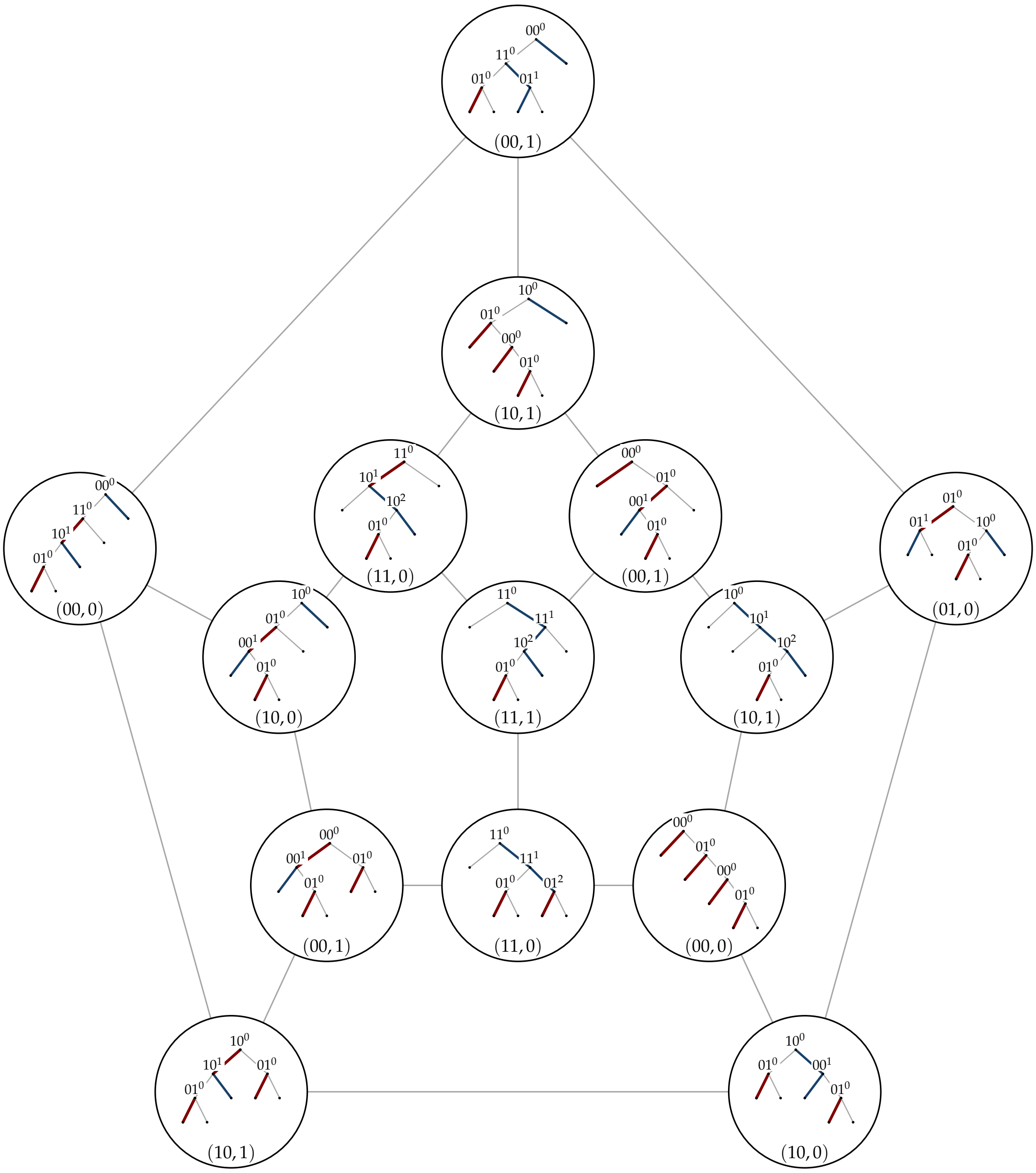}
\caption{The 8-colouring of $\mathcal A_6$. At each internal vertex, $ij^k$
records the label $(i,j)$ and path depth $k$; every leaf has label $(0,0)$.
Blue and red edges are preferred, and red edges are those counted
in $E(T)$. The pair below each tree is its colour $(\lab(T),d(T))$,
with $ij$ abbreviating $(i,j)$ and $d(T)$ the number of red edges
modulo $2$.}
\label{fig:associahedron6-colouring}
\end{figure}

For a vertex $v$ of $T$, the \defterm{path depth} $\defmath{\delta_T(v)}$
is the number of edges between $v$ and the top of the preferred path that
contains $v$. Let $\defmath{E(T)}$ be the set of preferred left edges $uv$ of
$T$ such that $u$ is the parent of $v$ and $\delta_T(u)$ is even. Define
\begin{equation*}
 \defmath{d(T)}\coloneqq |E(T)|\pmod2.
\end{equation*}
The statistic $d(T)$ adapts the path colour of Berendsohn, Cardinal,
Iacono and Kozma~\cite{BCIK2026}. They decompose a tree into heavy paths
and call a vertex left-heavy if its left subtree is larger than its
right subtree. Their path colour is the parity of the number of
left-heavy vertices whose distance from the top of their heavy path is
even. Here we use preferred paths in place of heavy paths and count
vertices whose preferred child is their left child.
Tan~\cite[Section~2]{Tan2026} uses the analogous statistic for
highlighted paths. Colour $T$ by
\begin{equation*}
 \defmath{\kappa(T)}\coloneqq \bigl(\lab(T),d(T)\bigr)\in\F^2\times\F.
\end{equation*}
The set $\F^2\times\F$ has $8$ elements. Thus our colouring uses four
labels and a single parity bit, whereas Tan's colouring uses three
labels and two parity bits.
Figure~\ref{fig:associahedron6-colouring} illustrates the computation of
these colours for $\mathcal A_6$.

We next show that the root label of a tree determines the label of any
rooted subtree, provided that the rest of the tree is fixed. When
replacing a rooted subtree at a vertex $z$, we identify $z$ with the root
of the replacement subtree, retaining its position as a left child,
right child or root of the whole tree. For a rotation, we call this
vertex the \defterm{local root}. Tan~\cite{Tan2026} uses the fact that, for a vertex and its
two children, any two of the three labels determine the third label. The
operation $Q$ has the same property, and Lemma~\ref{lem:context} is a
consequence of this property.

\begin{lemma}\label{lem:context}
Fix a tree and a vertex $z$, and let $X$ and $Y$ be trees. Replace the
rooted subtree at $z$ first by $X$ and then by $Y$. If the two resulting
trees have the same root label, then $\lab(X)=\lab(Y)$.
\end{lemma}
\begin{proof}
For all $a,c\in\F^2$, 
\begin{align*}
 Q(a,c)&=(a_0+a_1,\,a_1)+(c_0,\,c_0+c_1+1),\\
 Q(c,a)&=(a_0,\,a_0+a_1)+(c_0+c_1,\,c_1+1).
\end{align*}
The map $(a_0,a_1)\mapsto(a_0+a_1,a_1)$ is its own inverse, so this map
is a bijection of $\F^2$. Similarly, the map
$(a_0,a_1)\mapsto(a_0,a_0+a_1)$ is a bijection of $\F^2$. Hence
$Q(a,c)=Q(b,c)$ implies $a=b$. Likewise, $Q(c,a)=Q(c,b)$ implies $a=b$.

We proceed by induction on the depth of $z$ in the original
tree. If $z$ is the root, then the two resulting trees are $X$ and $Y$,
so the claim holds. Otherwise, let $p$ be the parent of $z$. By the induction
hypothesis, $p$ has the same label in both trees. The sibling subtree of
$z$ is unchanged, so its root also has the same label. As shown in the
first paragraph, the label of a parent and either child uniquely
determine the label of the other child. Hence $z$ has the same label in
both trees, giving $\lab(X)=\lab(Y)$.
\end{proof}

We now consider a rotation that preserves the label of its local root.

\begin{figure}[htbp]
\centering
\begin{tikzpicture}[x=1cm,y=1cm,
 vertex/.style={circle,fill=black,inner sep=1.7pt},
 preferred/.style={line width=1.3pt,draw=pathblue},
 other/.style={line width=.55pt,dashed,draw=black!60},
 continuation/.style={line width=.8pt,densely dotted,draw=pathblue},
 every node/.style={font=\small}]
 \coordinate (r) at (0,0);
 \coordinate (s) at (-.85,-.8);
 \coordinate (u) at (-1.65,-1.65);
 \coordinate (v) at (-.05,-1.65);
 \coordinate (w) at (1.25,-1.65);
 \draw[preferred] (r)--(s)--(v);
 \draw[other] (s)--(u) (r)--(w);
 \draw[continuation] (0,.55)--(r) (v)--(-.05,-2.18);
 \node[vertex,label=right:{$k$}] at (r) {};
 \node[vertex,label=left:{$k+1$}] at (s) {};
 \node[vertex,label=below left:{$U$}] at (u) {};
 \node[vertex,label=below left:{$V$},label=right:{$k+2$}] at (v) {};
 \node[vertex,label=below:{$W$}] at (w) {};
 \node[pathblue] at (-.55,-.18) {$L$};
 \node[pathblue] at (-.26,-1.10) {$R$};
 \draw[-{Stealth[length=2mm]},line width=.6pt] (1.65,-.72)--(2.45,-.72);
 \coordinate (rp) at (4.15,0);
 \coordinate (t) at (5,-.8);
 \coordinate (up) at (2.9,-1.65);
 \coordinate (vp) at (4.2,-1.65);
 \coordinate (wp) at (5.8,-1.65);
 \draw[preferred] (rp)--(t)--(vp);
 \draw[other] (rp)--(up) (t)--(wp);
 \draw[continuation] (4.15,.55)--(rp) (vp)--(4.2,-2.18);
 \node[vertex,label=left:{$k$}] at (rp) {};
 \node[vertex,label=right:{$k+1$}] at (t) {};
 \node[vertex,label=below:{$U$}] at (up) {};
 \node[vertex,label=below left:{$V$},label=right:{$k+2$}] at (vp) {};
 \node[vertex,label=below:{$W$}] at (wp) {};
 \node[pathblue] at (4.72,-.18) {$R$};
 \node[pathblue] at (4.42,-1.10) {$L$};
\end{tikzpicture}
\caption{A rotation that preserves the local root label. Solid edges are
preferred and dashed edges are unpreferred. Numbers are path depths,
where $k$ is the path depth of the local root. Dotted lines indicate
possible continuations of the preferred path.}
\label{fig:rotation}
\end{figure}
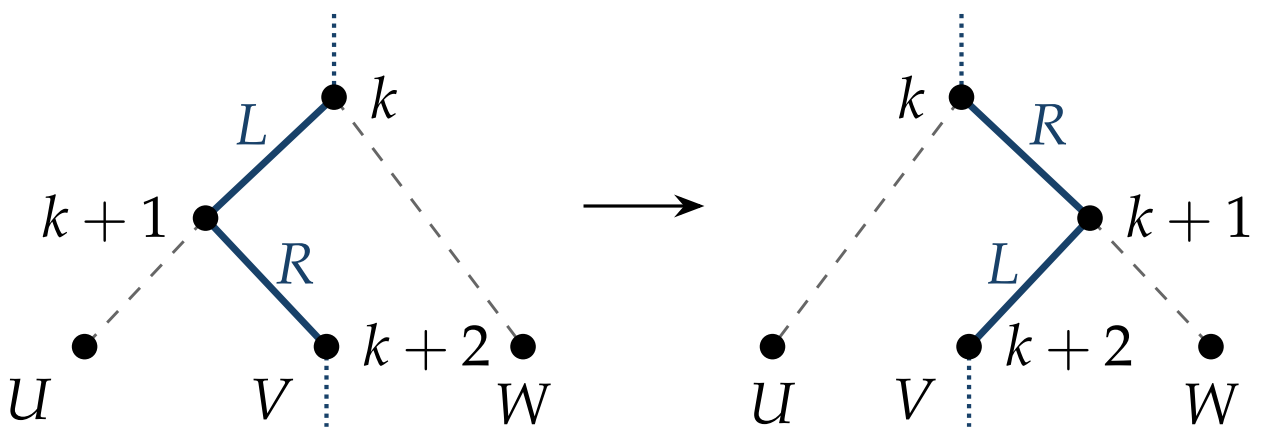

\Needspace{8\baselineskip}
\begin{lemma}\label{lem:swap}
Suppose that a rotation $((U,V),W)\to(U,(V,W))$ preserves the label of
its local root. Before the rotation, the path from the local root to the
root of $V$ consists of two preferred edges with word $LR$. After the
rotation, this path consists of two preferred edges with word $RL$.
Moreover, the edge that joins the root of $U$ to its parent is
unpreferred before and after the rotation. Similarly, the edge that joins
the root of $W$ to its parent is unpreferred before and after the
rotation.
\end{lemma}
\begin{proof}
Let $\defmath{a}\coloneqq\lab(U)$, $\defmath{b}\coloneqq\lab(V)$ and
$\defmath{c}\coloneqq\lab(W)$. Before the rotation, the left child of the
local root is the root of $(U,V)$, which has label
$\defmath{s}\coloneqq Q(a,b)$. After the rotation, the right child of the
local root is the root of $(V,W)$, which has label
$\defmath{t}\coloneqq Q(b,c)$. The labels of the local root before and
after the rotation are
\begin{align*}
 Q(Q(a,b),c)&=\bigl(a_0+b_1+c_0+1,\,
                  a_1+b_0+b_1+c_0+c_1\bigr),\\
 Q(a,Q(b,c))&=\bigl(a_0+a_1+b_0+b_1+c_0,\,
                  a_1+b_0+c_1\bigr).
\end{align*}
Adding these two expressions gives
\begin{equation*}
 Q(Q(a,b),c)+Q(a,Q(b,c))
       =\bigl(b_0+a_1+1,\ b_1+c_0\bigr).
\end{equation*}
Hence the label of the local root is the same before and after the
rotation if and only if
\begin{equation*}
       b_0=a_1+1\quad\text{and}\quad b_1=c_0.
\end{equation*}
Substituting into the definition of $D$ gives
\begin{align*}
 D(s,c)&=a_1+b_0+b_1+1+c_0=0,&
 D(a,b)&=a_1+b_0=1,\\
 D(a,t)&=a_1+b_0+b_1+c_0=1,&
 D(b,c)&=b_1+c_0=0.
\end{align*}
Before the rotation, $D(s,c)=0$, so the preferred child of the local root
is the root of $(U,V)$. Moreover, $D(a,b)=1$, so the preferred child of
the root of $(U,V)$ is the root of $V$. Thus the path from the local root
to the root of $V$ consists of two preferred edges with word $LR$. After
the rotation, $D(a,t)=1$, so the preferred child of the local root is the
root of $(V,W)$. Moreover, $D(b,c)=0$, so the preferred child of the root
of $(V,W)$ is the root of $V$. Thus the path from the local root to the
root of $V$ consists of two preferred edges with word $RL$. In both
trees, the preferred child of each of the two internal vertices on this
path lies on this path. Hence the edges that join the roots of $U$ and
$W$ to their parents are unpreferred in both trees.
Figure~\ref{fig:rotation} illustrates this change.
\end{proof}

\Needspace{8\baselineskip}
We can now prove that $\kappa$ is a colouring of $\An$.
\begin{proof}[Proof of Theorem~\ref{thm:main}]
Let $T$ and $T'$ be trees that differ by a rotation, and let $z$ be its
local root. After exchanging $T$ and $T'$ if necessary, the rotation
replaces the rooted subtree $((U,V),W)$ at $z$ in $T$ by $(U,(V,W))$ in
$T'$, leaving the rest of the tree unchanged. If
$\lab(T)\neq\lab(T')$, then $\kappa(T)\neq\kappa(T')$, since the first
coordinate of $\kappa$ is $\lab$. Now assume that $\lab(T)=\lab(T')$. By
Lemma~\ref{lem:context}, the local root $z$ has the same label in $T$
and $T'$.

Every vertex of $U$, $V$ and $W$ has the same label in $T$ and $T'$,
because these rooted subtrees are unchanged. Every vertex outside the
rooted subtree at $z$ also has the same label in $T$ and $T'$. This
statement holds because the label of such a vertex is computed from the
label of $z$ and from labels in unchanged rooted subtrees. Apart from
the local root and the internal child involved in the rotation, every
internal vertex retains its children and their labels. Its preferred
child therefore remains unchanged. Consequently, every preferred path that contains no vertex of
the rooted subtree at $z$ is the same in $T$ and $T'$. Moreover, the path
from the root of $T$ to $z$ is also the path from the root of $T'$ to
$z$. Every vertex on this path has the same path depth in $T$ and $T'$.
In particular, $z$ has the same path depth in $T$ and $T'$. Let
$\defmath{k}$ denote this path depth.

By Lemma~\ref{lem:swap}, the preferred path that contains $z$ continues
from $z$ to the root of $V$ with word $LR$ in $T$ and with word $RL$ in
$T'$. The part of this preferred path above $z$ is the same in $T$ and
$T'$. Every internal vertex of $V$ retains its preferred child. Hence the part of this
preferred path inside $V$ is the same in $T$ and $T'$. Since the root of
$V$ has path depth $k+2$ in both $T$ and $T'$, this part has the same
path depths in both trees. Every other preferred path inside $V$ is the
same in $T$ and $T'$. By Lemma~\ref{lem:swap}, the roots of $U$ and $W$
are tops of preferred paths in both $T$ and $T'$. Hence every preferred
path inside $U$ or $W$ is the same in $T$ and $T'$. Moreover, every
vertex of $U$ and $W$ has the same path depth in $T$ and $T'$.

We now compare the preferred left edges of $T$ and $T'$. By
Lemma~\ref{lem:swap}, the edges that join the roots of $U$ and $W$ to
their parents are unpreferred in both $T$ and $T'$. Every other edge
outside the two-edge segment from $z$ to the root of $V$ retains its
preferred status and the path depth of its upper endpoint. Thus all
contributions to the parity bit outside this segment remain unchanged.

The path from $z$ to the root of $V$ contains exactly one left edge in $T$
and exactly one left edge in $T'$. In $T$, the upper endpoint of
this left edge is $z$, which has path depth $k$. In $T'$, the upper
endpoint of this left edge is the root of $(V,W)$, which has path depth
$k+1$. Since $k$ and $k+1$ have opposite parity, exactly one of these two
left edges contributes to the corresponding count. Therefore
$d(T')=d(T)+1$ in $\F$, and so $\kappa(T)\neq\kappa(T')$. Thus $\kappa$
is a colouring of $\An$ with at most $8$ colours.
\end{proof}

We finish this section by describing two variations of the colouring
$\kappa$.
\begin{itemize}
\item For each $m\geqslant1$, fix a sequence of $m$ labels in $\F^2$. In
every tree with $m$ leaves, assign these labels to the leaves from left
to right. Then label the internal vertices using $Q$ as before. With this
labelling, the label of a vertex depends on its rooted subtree and also
on the positions of the leaves of that rooted subtree.
The following analogue of Lemma~\ref{lem:context} holds. Fix a tree and
a vertex $z$, and let $X$ and $Y$ have the same number of leaves.
Replace the rooted subtree at $z$ first by $X$ and then by $Y$. If the
two resulting trees have the same root label, then $z$ has the same
label in both trees. The same proof applies, because every leaf
outside the rooted subtree at $z$ keeps its position and its label. Now
let $T$ and $T'$ be two trees that differ by a rotation. A rotation
preserves the number of leaves and their left-to-right order. Hence
every leaf has the same label in $T$ and $T'$. Therefore the proof of
Theorem~\ref{thm:main} applies to the resulting colouring.
\item An \defterm{inversion} in a word is a pair of positions $i<j$ such
that the letter at position $i$ is $R$ and the letter at position $j$ is
$L$. The \defterm{inversion parity} of a collection of words is the total
number of inversions in these words modulo $2$. Oum and
Wood~\cite[Section~3]{OW2026} use the inversion parity of the words of
their paths. Tan~\cite[Section~2]{Tan2026} observes that this inversion
parity can replace the analogue of the path colour in Tan's colouring. The same
replacement works for preferred paths. Indeed, let $\defmath{d'(T)}$ be
the inversion parity of the words of the preferred paths of $T$.
Replacing two adjacent letters $LR$ by $RL$ changes the number of
inversions in a word by exactly $1$. Now let $T$ and $T'$ be two trees
that differ by a rotation that preserves the label of the root. By the
proof of Theorem~\ref{thm:main}, every preferred path other than the one
that contains the local root has the same word in $T$ and $T'$. The word
of the preferred path that contains the local root changes by replacing
two adjacent letters $LR$ by $RL$. Hence $d'(T')\neq d'(T)$, and the proof of
Theorem~\ref{thm:main} shows that $T\mapsto(\lab(T),d'(T))$ is also a
colouring of $\An$.
\end{itemize}

\section{A Lower Bound for Label--Parity Colourings}\label{sec:lower}

The colouring $\kappa$ used in Section~\ref{sec:proof} pairs a recursive label with one extra bit. We
call a colouring of this form a \defterm{label--parity colouring}. This
section considers label--parity colourings in which every leaf has the
same label and the labels come from a quasigroup. For such a colouring,
we show that the set of possible colours has at least $8$ elements.

A \defterm{quasigroup} is a set $S$ together with a binary operation
$\cdot$ on $S$ such that, for all $a,b\in S$, each of the equations
$a\cdot y=b$ and $y\cdot a=b$ has exactly one solution $y\in S$.
Equivalently, every row and every column of the operation table of $S$
is a permutation of $S$. A finite table with this property is a
\defterm{Latin square}. The proof of Lemma~\ref{lem:context} shows that
every row and every column of the operation table of $Q$ is a bijection
of $\F^2$. Hence $(\F^2,Q)$ is a quasigroup. Fix a finite quasigroup $(S,\cdot)$ and an
element $x\in S$. Every leaf receives the label $x$. Consider an internal
vertex whose left child has label $a$ and whose right child has label
$b$. This internal vertex receives the label $a\cdot b$. As before,
$\lab(T)$ denotes the label of the root of a tree $T$.

Let $\defmath{S'}$ be the set of labels $\lab(T)$ over all trees $T$.
The set $S'$ is closed under $\cdot$: if $a,b\in S'$ are realised by
trees $U,V$, respectively, then $a\cdot b=\lab((U,V))\in S'$.
Each row and each column of the operation table
restricted to $S'$ is injective and hence, since $S'$ is finite, is a
bijection of $S'$. Thus $S'$ is a quasigroup.

A triple $(a,b,c)\in S'^3$ is \defterm{associative} if
$(a\cdot b)\cdot c=a\cdot(b\cdot c)$. The extra parity bit in our proof
is necessary, since if the root label alone gave a colouring of $\An$
for every $n$, then
$S'$ could contain no associative triple. Indeed, trees $A$, $B$ and $C$
realising such a triple would give adjacent trees $((A,B),C)$ and
$(A,(B,C))$ with the same root label. However, as noted by Dr\'apal and
Wanless~\cite{DW2021}, every finite quasigroup has at least as many
associative triples as elements. Thus the root label alone cannot colour
$\An$ for every $n$. We next show that combining the root label with one
extra bit requires at least four realised labels.

\begin{theorem}\label{thm:lower}
Let $(S,\cdot)$, $x$ and $\lab$ be as above. Let $\defmath{\beta}$ be an
arbitrary function from the set of ordered full binary trees to $\F$.
Suppose that, for every $n\geqslant3$,
the map $T\mapsto(\lab(T),\beta(T))$ is a colouring of $\An$. Then
$|S'|\geqslant4$. Hence the map $T\mapsto(\lab(T),\beta(T))$ takes values
in the set $S'\times\F$, which has at least $8$ elements.
\end{theorem}

\begin{proof}
Assume, for a contradiction, that $|S'|\leqslant3$.

Consider a rotation
$((U,V),W)\to(U,(V,W))$ inside a tree, and suppose that
$(\lab(U),\lab(V),\lab(W))$ is associative. Then the label of the local
root is the same before and after the rotation. The rest of the tree is
unchanged. Hence the root label of the whole tree is the same before and
after the rotation. Since $T\mapsto(\lab(T),\beta(T))$ is a colouring and
$\beta$ takes values in $\{0,1\}$, $\beta$ must flip at every such
rotation. An odd cycle of such rotations would force $\beta$ to return
to its starting vertex with the opposite value, a contradiction. It
therefore suffices to find such an odd cycle in some $\An$.

Let $A$, $B$, $C$ and $E$ be trees with labels $a$, $b$, $c$ and $e$,
respectively. The five trees in Figure~\ref{fig:pentagon}, read in the
numbered cyclic order, form a cycle of length $5$. We call this
cycle the \defterm{pentagon}. The five rotations of the pentagon occur at the label
triples
\begin{equation}\label{eq:pentagon}
 (a,b,c),\quad (a,b\cdot c,e),\quad (b,c,e),\quad (a,b,c\cdot e),\quad
 (a\cdot b,c,e).
\end{equation}
If all five triples in~\eqref{eq:pentagon} are associative, then the
pentagon is the required odd cycle.

\begin{figure}[htbp]
\centering
\begin{tikzpicture}[x=1cm,y=1cm,scale=.9,line cap=round,line join=round,
 every node/.style={font=\small}]
\coordinate (a0) at (0.00000,2.80000);
\node at (0.00000,4.33000) {$ 1 $};
\coordinate (a1) at (3.32870,0.86525);
\node at (4.78381,1.33804) {$ 2 $};
\coordinate (a2) at (2.05725,-2.26525);
\node at (2.95656,-3.50304) {$ 3 $};
\coordinate (a3) at (-2.05725,-2.26525);
\node at (-2.95656,-3.50304) {$ 4 $};
\coordinate (a4) at (-3.32870,0.86525);
\node at (-4.78381,1.33804) {$ 5 $};
\draw[pathblue!65,line width=.65pt,-{Stealth[length=1.6mm]},shorten <=1.17cm,shorten >=1.17cm] (a0)--(a1);
\draw[pathblue!65,line width=.65pt,-{Stealth[length=1.6mm]},shorten <=1.17cm,shorten >=1.17cm] (a1)--(a2);
\draw[pathblue!65,line width=.65pt,-{Stealth[length=1.6mm]},shorten <=1.17cm,shorten >=1.17cm] (a2)--(a3);
\draw[pathblue!65,line width=.65pt,-{Stealth[length=1.6mm]},shorten <=1.17cm,shorten >=1.17cm] (a3)--(a4);
\draw[pathblue!65,line width=.65pt,-{Stealth[length=1.6mm]},shorten <=1.17cm,shorten >=1.17cm] (a4)--(a0);
\begin{scope}[shift={(a0)}]
\draw[draw=black,fill=white,line width=.55pt] (0,0) circle[radius=1.25];
\coordinate (t0) at (0.32500,0.81000);
\coordinate (t1) at (-0.13000,0.39000);
\coordinate (t2) at (-0.52000,-0.03000);
\coordinate (t3) at (-0.78000,-0.45000);
\coordinate (t4) at (-0.26000,-0.45000);
\coordinate (t5) at (0.26000,-0.03000);
\coordinate (t6) at (0.78000,0.39000);
\draw[black,line width=.5pt] (t2)--(t3);
\draw[black,line width=.5pt] (t2)--(t4);
\draw[black,line width=.5pt] (t1)--(t2);
\draw[black,line width=.5pt] (t1)--(t5);
\draw[black,line width=.5pt] (t0)--(t1);
\draw[black,line width=.5pt] (t0)--(t6);
\fill (t0) circle[radius=.022];
\fill (t1) circle[radius=.022];
\fill (t2) circle[radius=.022];
\draw[draw=black,fill=black!7,line width=.45pt] (t3)--++(-.20,-.36)--++(.40,0)--cycle;
\node[font=\tiny,inner sep=0pt] at ([yshift=-2.5mm]t3) {$A$};
\draw[draw=black,fill=black!7,line width=.45pt] (t4)--++(-.20,-.36)--++(.40,0)--cycle;
\node[font=\tiny,inner sep=0pt] at ([yshift=-2.5mm]t4) {$B$};
\draw[draw=black,fill=black!7,line width=.45pt] (t5)--++(-.20,-.36)--++(.40,0)--cycle;
\node[font=\tiny,inner sep=0pt] at ([yshift=-2.5mm]t5) {$C$};
\draw[draw=black,fill=black!7,line width=.45pt] (t6)--++(-.20,-.36)--++(.40,0)--cycle;
\node[font=\tiny,inner sep=0pt] at ([yshift=-2.5mm]t6) {$E$};
\end{scope}
\begin{scope}[shift={(a1)}]
\draw[draw=black,fill=white,line width=.55pt] (0,0) circle[radius=1.25];
\coordinate (t0) at (0.19500,0.81000);
\coordinate (t1) at (-0.39000,0.39000);
\coordinate (t2) at (-0.78000,-0.03000);
\coordinate (t3) at (0.00000,-0.03000);
\coordinate (t4) at (-0.26000,-0.45000);
\coordinate (t5) at (0.26000,-0.45000);
\coordinate (t6) at (0.78000,0.39000);
\draw[black,line width=.5pt] (t3)--(t4);
\draw[black,line width=.5pt] (t3)--(t5);
\draw[black,line width=.5pt] (t1)--(t2);
\draw[black,line width=.5pt] (t1)--(t3);
\draw[black,line width=.5pt] (t0)--(t1);
\draw[black,line width=.5pt] (t0)--(t6);
\fill (t0) circle[radius=.022];
\fill (t1) circle[radius=.022];
\draw[draw=black,fill=black!7,line width=.45pt] (t2)--++(-.20,-.36)--++(.40,0)--cycle;
\node[font=\tiny,inner sep=0pt] at ([yshift=-2.5mm]t2) {$A$};
\fill (t3) circle[radius=.022];
\draw[draw=black,fill=black!7,line width=.45pt] (t4)--++(-.20,-.36)--++(.40,0)--cycle;
\node[font=\tiny,inner sep=0pt] at ([yshift=-2.5mm]t4) {$B$};
\draw[draw=black,fill=black!7,line width=.45pt] (t5)--++(-.20,-.36)--++(.40,0)--cycle;
\node[font=\tiny,inner sep=0pt] at ([yshift=-2.5mm]t5) {$C$};
\draw[draw=black,fill=black!7,line width=.45pt] (t6)--++(-.20,-.36)--++(.40,0)--cycle;
\node[font=\tiny,inner sep=0pt] at ([yshift=-2.5mm]t6) {$E$};
\end{scope}
\begin{scope}[shift={(a2)}]
\draw[draw=black,fill=white,line width=.55pt] (0,0) circle[radius=1.25];
\coordinate (t0) at (-0.19500,0.81000);
\coordinate (t1) at (-0.78000,0.39000);
\coordinate (t2) at (0.39000,0.39000);
\coordinate (t3) at (0.00000,-0.03000);
\coordinate (t4) at (-0.26000,-0.45000);
\coordinate (t5) at (0.26000,-0.45000);
\coordinate (t6) at (0.78000,-0.03000);
\draw[black,line width=.5pt] (t3)--(t4);
\draw[black,line width=.5pt] (t3)--(t5);
\draw[black,line width=.5pt] (t2)--(t3);
\draw[black,line width=.5pt] (t2)--(t6);
\draw[black,line width=.5pt] (t0)--(t1);
\draw[black,line width=.5pt] (t0)--(t2);
\fill (t0) circle[radius=.022];
\draw[draw=black,fill=black!7,line width=.45pt] (t1)--++(-.20,-.36)--++(.40,0)--cycle;
\node[font=\tiny,inner sep=0pt] at ([yshift=-2.5mm]t1) {$A$};
\fill (t2) circle[radius=.022];
\fill (t3) circle[radius=.022];
\draw[draw=black,fill=black!7,line width=.45pt] (t4)--++(-.20,-.36)--++(.40,0)--cycle;
\node[font=\tiny,inner sep=0pt] at ([yshift=-2.5mm]t4) {$B$};
\draw[draw=black,fill=black!7,line width=.45pt] (t5)--++(-.20,-.36)--++(.40,0)--cycle;
\node[font=\tiny,inner sep=0pt] at ([yshift=-2.5mm]t5) {$C$};
\draw[draw=black,fill=black!7,line width=.45pt] (t6)--++(-.20,-.36)--++(.40,0)--cycle;
\node[font=\tiny,inner sep=0pt] at ([yshift=-2.5mm]t6) {$E$};
\end{scope}
\begin{scope}[shift={(a3)}]
\draw[draw=black,fill=white,line width=.55pt] (0,0) circle[radius=1.25];
\coordinate (t0) at (-0.32500,0.81000);
\coordinate (t1) at (-0.78000,0.39000);
\coordinate (t2) at (0.13000,0.39000);
\coordinate (t3) at (-0.26000,-0.03000);
\coordinate (t4) at (0.52000,-0.03000);
\coordinate (t5) at (0.26000,-0.45000);
\coordinate (t6) at (0.78000,-0.45000);
\draw[black,line width=.5pt] (t4)--(t5);
\draw[black,line width=.5pt] (t4)--(t6);
\draw[black,line width=.5pt] (t2)--(t3);
\draw[black,line width=.5pt] (t2)--(t4);
\draw[black,line width=.5pt] (t0)--(t1);
\draw[black,line width=.5pt] (t0)--(t2);
\fill (t0) circle[radius=.022];
\draw[draw=black,fill=black!7,line width=.45pt] (t1)--++(-.20,-.36)--++(.40,0)--cycle;
\node[font=\tiny,inner sep=0pt] at ([yshift=-2.5mm]t1) {$A$};
\fill (t2) circle[radius=.022];
\draw[draw=black,fill=black!7,line width=.45pt] (t3)--++(-.20,-.36)--++(.40,0)--cycle;
\node[font=\tiny,inner sep=0pt] at ([yshift=-2.5mm]t3) {$B$};
\fill (t4) circle[radius=.022];
\draw[draw=black,fill=black!7,line width=.45pt] (t5)--++(-.20,-.36)--++(.40,0)--cycle;
\node[font=\tiny,inner sep=0pt] at ([yshift=-2.5mm]t5) {$C$};
\draw[draw=black,fill=black!7,line width=.45pt] (t6)--++(-.20,-.36)--++(.40,0)--cycle;
\node[font=\tiny,inner sep=0pt] at ([yshift=-2.5mm]t6) {$E$};
\end{scope}
\begin{scope}[shift={(a4)}]
\draw[draw=black,fill=white,line width=.55pt] (0,0) circle[radius=1.25];
\coordinate (t0) at (0.00000,0.60000);
\coordinate (t1) at (-0.52000,0.18000);
\coordinate (t2) at (-0.78000,-0.24000);
\coordinate (t3) at (-0.26000,-0.24000);
\coordinate (t4) at (0.52000,0.18000);
\coordinate (t5) at (0.26000,-0.24000);
\coordinate (t6) at (0.78000,-0.24000);
\draw[black,line width=.5pt] (t1)--(t2);
\draw[black,line width=.5pt] (t1)--(t3);
\draw[black,line width=.5pt] (t4)--(t5);
\draw[black,line width=.5pt] (t4)--(t6);
\draw[black,line width=.5pt] (t0)--(t1);
\draw[black,line width=.5pt] (t0)--(t4);
\fill (t0) circle[radius=.022];
\fill (t1) circle[radius=.022];
\draw[draw=black,fill=black!7,line width=.45pt] (t2)--++(-.20,-.36)--++(.40,0)--cycle;
\node[font=\tiny,inner sep=0pt] at ([yshift=-2.5mm]t2) {$A$};
\draw[draw=black,fill=black!7,line width=.45pt] (t3)--++(-.20,-.36)--++(.40,0)--cycle;
\node[font=\tiny,inner sep=0pt] at ([yshift=-2.5mm]t3) {$B$};
\fill (t4) circle[radius=.022];
\draw[draw=black,fill=black!7,line width=.45pt] (t5)--++(-.20,-.36)--++(.40,0)--cycle;
\node[font=\tiny,inner sep=0pt] at ([yshift=-2.5mm]t5) {$C$};
\draw[draw=black,fill=black!7,line width=.45pt] (t6)--++(-.20,-.36)--++(.40,0)--cycle;
\node[font=\tiny,inner sep=0pt] at ([yshift=-2.5mm]t6) {$E$};
\end{scope}
\end{tikzpicture}
\caption{The pentagon. Each
arrow is a rotation, and the shaded triangles represent the unchanged
rooted subtrees $A$, $B$, $C$ and $E$.}
\label{fig:pentagon}
\end{figure}
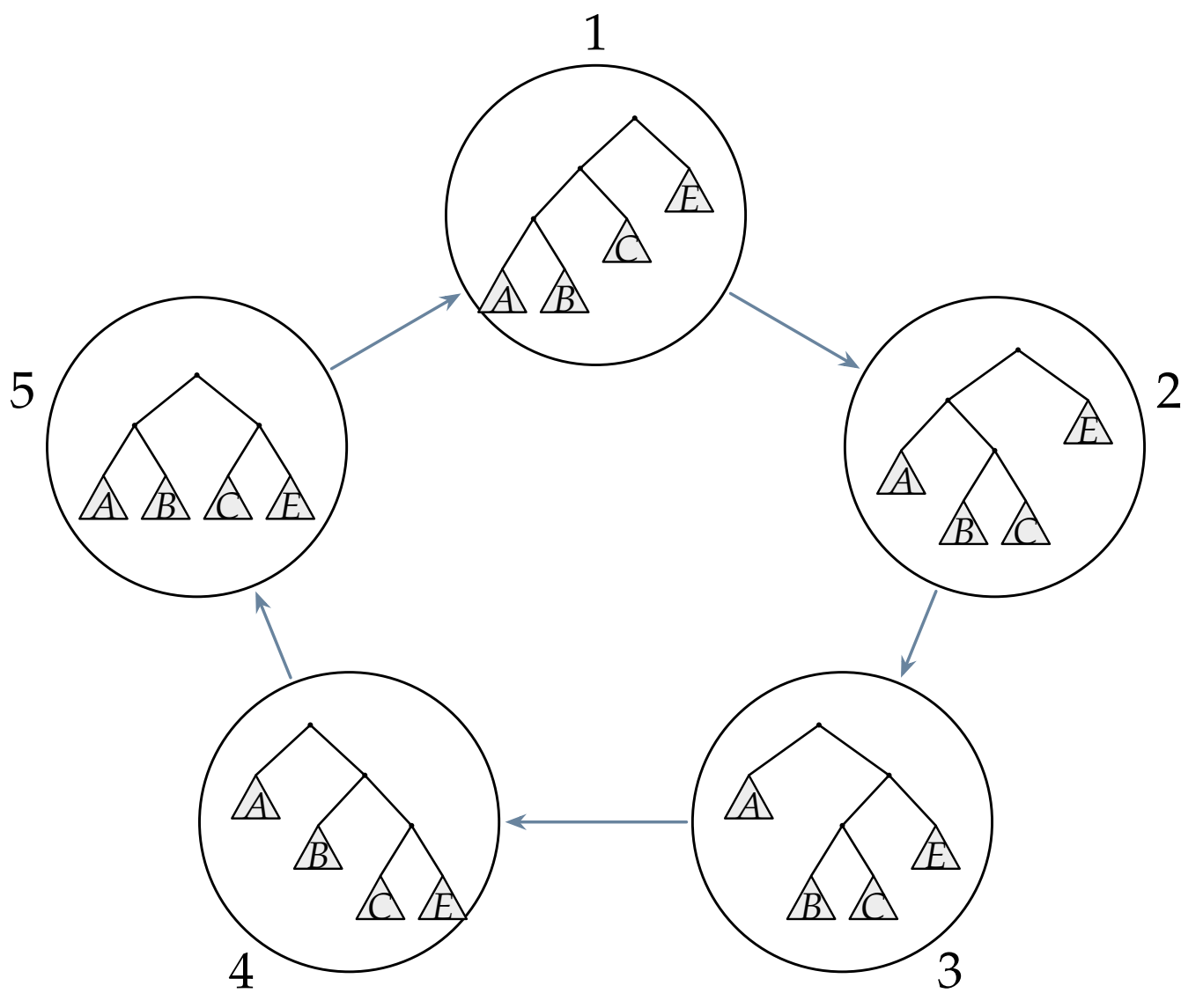

\emph{Case.}~$|S'|\leqslant2$: If $|S'|=1$, then every triple is
associative. If $|S'|=2$, then we identify $S'$ with $\F$. There
are only two Latin squares of order $2$, namely the tables of the
operations $a\cdot b=a+b+w$ with $w\in\F$. Each of these
operations is associative. Hence, in both cases, every triple is
associative. Thus the pentagon gives a contradiction.

\emph{Case.}~$|S'|=3$: We identify $S'$ with $\mathbb Z_3$. Since every
element of $S'$ is the label of some tree, the trees $A$, $B$, $C$ and
$E$ in the pentagon can have any prescribed labels in $\mathbb Z_3$.
There are at most $12$ Latin squares of order $3$. Indeed, the first row
is one of $6$ permutations. The second row is one of the $2$ permutations
that disagree with the first row in every column. The third row is then
determined. Now consider the $12$ operations $a\cdot b=ua+vb+w$ with
$u,v\in\{1,-1\}$ and $w\in\mathbb Z_3$. The table of each of these
operations is a Latin square, because $u$ and $v$ are invertible in
$\mathbb Z_3$. These $12$ tables are distinct, because $w=0\cdot0$,
$u=1\cdot0-w$ and $v=0\cdot1-w$. Hence every Latin square of order $3$ is
the table of one of these $12$ operations. In particular, the operation
of $S'$ is $a\cdot b=ua+vb+w$ for some $u,v\in\{1,-1\}$ and
$w\in\mathbb Z_3$. Direct expansion gives
\begin{equation*}
 (a\cdot b)\cdot c-a\cdot(b\cdot c)=(1-u)a+(v-1)c+(u-v)w.
\end{equation*}
In particular, whether $(a,b,c)$ is associative does not depend on $b$.
We consider four subcases, according to the values of $u$ and $v$.

\emph{Subcase.}~$u=v=1$: Every triple is associative. Hence the pentagon
gives a contradiction.

\emph{Subcase.}~$u=1$ and $v=-1$: The triple $(a,b,c)$ is associative if and
only if $c=w$. Choose $c=e=w$. Then $c\cdot e=w-w+w=w$. Hence every triple
in~\eqref{eq:pentagon} has third entry $w$. Thus every triple
in~\eqref{eq:pentagon} is associative, and the pentagon gives a
contradiction.

\emph{Subcase.}~$u=-1$ and $v=1$: The triple $(a,b,c)$ is associative if and
only if $a=w$. Choose $a=b=w$. Then $a\cdot b=-w+w+w=w$. Hence every
triple in~\eqref{eq:pentagon} has first entry $w$. Thus every triple
in~\eqref{eq:pentagon} is associative, and the pentagon gives a
contradiction.

\emph{Subcase.}~$u=v=-1$: The triple $(a,b,c)$ is associative if and only if
$a=c$. The pentagon gives a contradiction only when $w=0$. Hence, in this
case, we use a cycle of length $9$ instead of the pentagon. Let $\defmath{t}$ denote the tree with one leaf, and
abbreviate $(U,V)$ to $UV$. Then $\lab(t)=x$ and
$\lab(tt)=x+w$. Moreover,
$\lab(t(tt))=\lab((tt)t)=x$. The nine trees with $5$ leaves in
Figure~\ref{fig:nine-walk}, read in the numbered cyclic order, form a
cycle of length $9$ in $\mathcal A_6$.

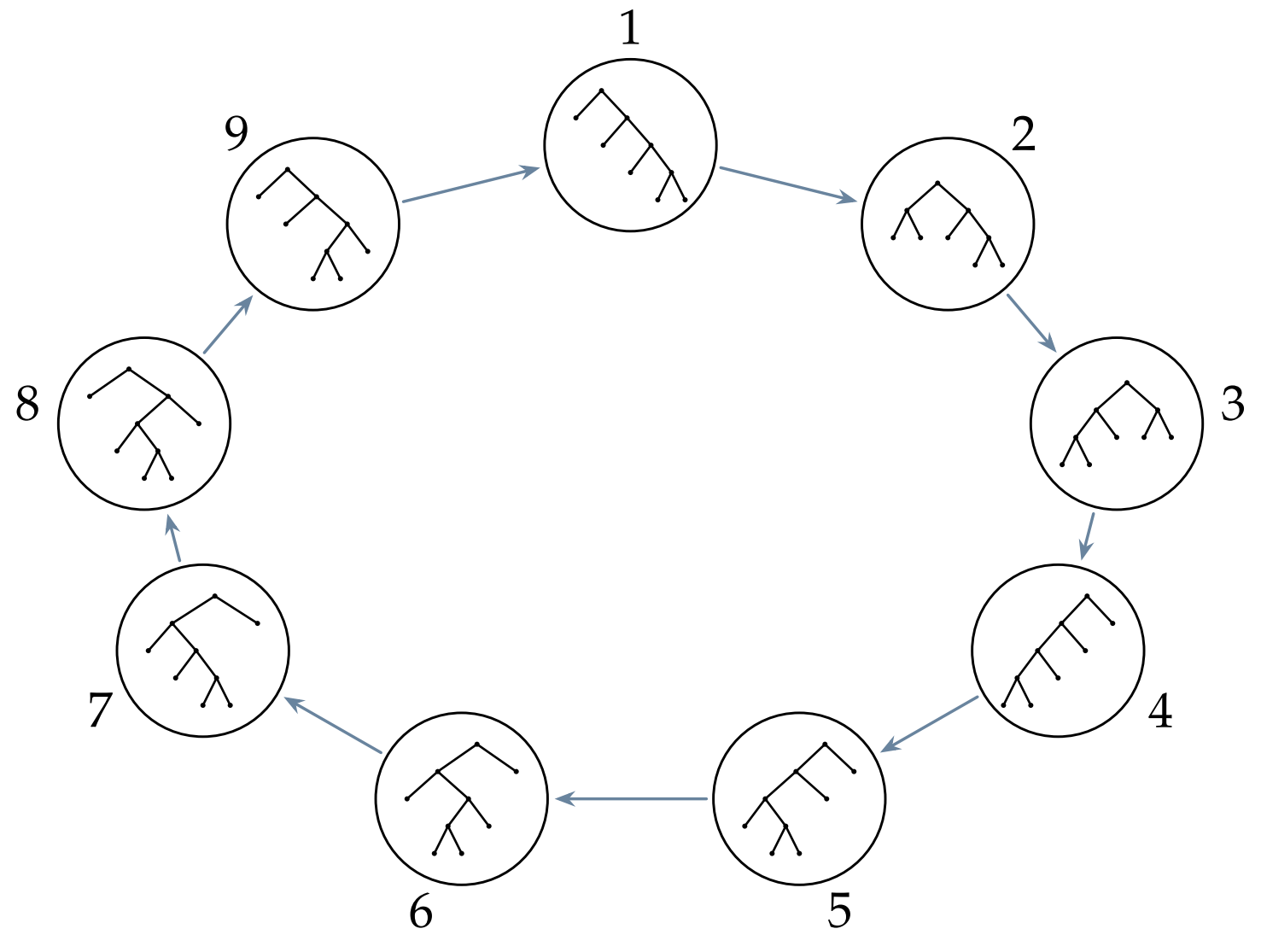
\begin{figure}[htbp]
\centering
\begin{tikzpicture}[x=1cm,y=1cm,scale=.9,line cap=round,line join=round,
 every node/.style={font=\small}]
\coordinate (a0) at (0.00000,2.90000);
\node at (0.00000,3.92000) {$ 1 $};
\coordinate (a1) at (2.73185,2.22153);
\node at (3.38749,3.00289) {$ 2 $};
\coordinate (a2) at (4.18543,0.50358);
\node at (5.18994,0.68070) {$ 3 $};
\coordinate (a3) at (3.68061,-1.45000);
\node at (4.56395,-1.96000) {$ 4 $};
\coordinate (a4) at (1.45359,-2.72511);
\node at (1.80245,-3.68360) {$ 5 $};
\coordinate (a5) at (-1.45359,-2.72511);
\node at (-1.80245,-3.68360) {$ 6 $};
\coordinate (a6) at (-3.68061,-1.45000);
\node at (-4.56395,-1.96000) {$ 7 $};
\coordinate (a7) at (-4.18543,0.50358);
\node at (-5.18994,0.68070) {$ 8 $};
\coordinate (a8) at (-2.73185,2.22153);
\node at (-3.38749,3.00289) {$ 9 $};
\draw[pathblue!65,line width=.65pt,-{Stealth[length=1.6mm]},shorten <=.72cm,shorten >=.72cm] (a0)--(a1);
\draw[pathblue!65,line width=.65pt,-{Stealth[length=1.6mm]},shorten <=.72cm,shorten >=.72cm] (a1)--(a2);
\draw[pathblue!65,line width=.65pt,-{Stealth[length=1.6mm]},shorten <=.72cm,shorten >=.72cm] (a2)--(a3);
\draw[pathblue!65,line width=.65pt,-{Stealth[length=1.6mm]},shorten <=.72cm,shorten >=.72cm] (a3)--(a4);
\draw[pathblue!65,line width=.65pt,-{Stealth[length=1.6mm]},shorten <=.72cm,shorten >=.72cm] (a4)--(a5);
\draw[pathblue!65,line width=.65pt,-{Stealth[length=1.6mm]},shorten <=.72cm,shorten >=.72cm] (a5)--(a6);
\draw[pathblue!65,line width=.65pt,-{Stealth[length=1.6mm]},shorten <=.72cm,shorten >=.72cm] (a6)--(a7);
\draw[pathblue!65,line width=.65pt,-{Stealth[length=1.6mm]},shorten <=.72cm,shorten >=.72cm] (a7)--(a8);
\draw[pathblue!65,line width=.65pt,-{Stealth[length=1.6mm]},shorten <=.72cm,shorten >=.72cm] (a8)--(a0);
\begin{scope}[shift={(a0)}]
\draw[draw=black,fill=white,line width=.55pt] (0,0) circle[radius=.74];
\coordinate (t0) at (-0.24969,0.47000);
\coordinate (t1) at (-0.47000,0.23500);
\coordinate (t2) at (-0.02937,0.23500);
\coordinate (t3) at (-0.23500,0.00000);
\coordinate (t4) at (0.17625,0.00000);
\coordinate (t5) at (0.00000,-0.23500);
\coordinate (t6) at (0.35250,-0.23500);
\coordinate (t7) at (0.23500,-0.47000);
\coordinate (t8) at (0.47000,-0.47000);
\draw[black,line width=.5pt] (t6)--(t7);
\draw[black,line width=.5pt] (t6)--(t8);
\draw[black,line width=.5pt] (t4)--(t5);
\draw[black,line width=.5pt] (t4)--(t6);
\draw[black,line width=.5pt] (t2)--(t3);
\draw[black,line width=.5pt] (t2)--(t4);
\draw[black,line width=.5pt] (t0)--(t1);
\draw[black,line width=.5pt] (t0)--(t2);
\foreach \v in {0,...,8}{\fill (t\v) circle[radius=.022];}
\end{scope}
\begin{scope}[shift={(a1)}]
\draw[draw=black,fill=white,line width=.55pt] (0,0) circle[radius=.74];
\coordinate (t0) at (-0.08812,0.35250);
\coordinate (t1) at (-0.35250,0.11750);
\coordinate (t2) at (-0.47000,-0.11750);
\coordinate (t3) at (-0.23500,-0.11750);
\coordinate (t4) at (0.17625,0.11750);
\coordinate (t5) at (0.00000,-0.11750);
\coordinate (t6) at (0.35250,-0.11750);
\coordinate (t7) at (0.23500,-0.35250);
\coordinate (t8) at (0.47000,-0.35250);
\draw[black,line width=.5pt] (t1)--(t2);
\draw[black,line width=.5pt] (t1)--(t3);
\draw[black,line width=.5pt] (t6)--(t7);
\draw[black,line width=.5pt] (t6)--(t8);
\draw[black,line width=.5pt] (t4)--(t5);
\draw[black,line width=.5pt] (t4)--(t6);
\draw[black,line width=.5pt] (t0)--(t1);
\draw[black,line width=.5pt] (t0)--(t4);
\foreach \v in {0,...,8}{\fill (t\v) circle[radius=.022];}
\end{scope}
\begin{scope}[shift={(a2)}]
\draw[draw=black,fill=white,line width=.55pt] (0,0) circle[radius=.74];
\coordinate (t0) at (0.08812,0.35250);
\coordinate (t1) at (-0.17625,0.11750);
\coordinate (t2) at (-0.35250,-0.11750);
\coordinate (t3) at (-0.47000,-0.35250);
\coordinate (t4) at (-0.23500,-0.35250);
\coordinate (t5) at (0.00000,-0.11750);
\coordinate (t6) at (0.35250,0.11750);
\coordinate (t7) at (0.23500,-0.11750);
\coordinate (t8) at (0.47000,-0.11750);
\draw[black,line width=.5pt] (t2)--(t3);
\draw[black,line width=.5pt] (t2)--(t4);
\draw[black,line width=.5pt] (t1)--(t2);
\draw[black,line width=.5pt] (t1)--(t5);
\draw[black,line width=.5pt] (t6)--(t7);
\draw[black,line width=.5pt] (t6)--(t8);
\draw[black,line width=.5pt] (t0)--(t1);
\draw[black,line width=.5pt] (t0)--(t6);
\foreach \v in {0,...,8}{\fill (t\v) circle[radius=.022];}
\end{scope}
\begin{scope}[shift={(a3)}]
\draw[draw=black,fill=white,line width=.55pt] (0,0) circle[radius=.74];
\coordinate (t0) at (0.24969,0.47000);
\coordinate (t1) at (0.02937,0.23500);
\coordinate (t2) at (-0.17625,0.00000);
\coordinate (t3) at (-0.35250,-0.23500);
\coordinate (t4) at (-0.47000,-0.47000);
\coordinate (t5) at (-0.23500,-0.47000);
\coordinate (t6) at (0.00000,-0.23500);
\coordinate (t7) at (0.23500,0.00000);
\coordinate (t8) at (0.47000,0.23500);
\draw[black,line width=.5pt] (t3)--(t4);
\draw[black,line width=.5pt] (t3)--(t5);
\draw[black,line width=.5pt] (t2)--(t3);
\draw[black,line width=.5pt] (t2)--(t6);
\draw[black,line width=.5pt] (t1)--(t2);
\draw[black,line width=.5pt] (t1)--(t7);
\draw[black,line width=.5pt] (t0)--(t1);
\draw[black,line width=.5pt] (t0)--(t8);
\foreach \v in {0,...,8}{\fill (t\v) circle[radius=.022];}
\end{scope}
\begin{scope}[shift={(a4)}]
\draw[draw=black,fill=white,line width=.55pt] (0,0) circle[radius=.74];
\coordinate (t0) at (0.22031,0.47000);
\coordinate (t1) at (-0.02937,0.23500);
\coordinate (t2) at (-0.29375,0.00000);
\coordinate (t3) at (-0.47000,-0.23500);
\coordinate (t4) at (-0.11750,-0.23500);
\coordinate (t5) at (-0.23500,-0.47000);
\coordinate (t6) at (0.00000,-0.47000);
\coordinate (t7) at (0.23500,0.00000);
\coordinate (t8) at (0.47000,0.23500);
\draw[black,line width=.5pt] (t4)--(t5);
\draw[black,line width=.5pt] (t4)--(t6);
\draw[black,line width=.5pt] (t2)--(t3);
\draw[black,line width=.5pt] (t2)--(t4);
\draw[black,line width=.5pt] (t1)--(t2);
\draw[black,line width=.5pt] (t1)--(t7);
\draw[black,line width=.5pt] (t0)--(t1);
\draw[black,line width=.5pt] (t0)--(t8);
\foreach \v in {0,...,8}{\fill (t\v) circle[radius=.022];}
\end{scope}
\begin{scope}[shift={(a5)}]
\draw[draw=black,fill=white,line width=.55pt] (0,0) circle[radius=.74];
\coordinate (t0) at (0.13219,0.47000);
\coordinate (t1) at (-0.20563,0.23500);
\coordinate (t2) at (-0.47000,0.00000);
\coordinate (t3) at (0.05875,0.00000);
\coordinate (t4) at (-0.11750,-0.23500);
\coordinate (t5) at (-0.23500,-0.47000);
\coordinate (t6) at (0.00000,-0.47000);
\coordinate (t7) at (0.23500,-0.23500);
\coordinate (t8) at (0.47000,0.23500);
\draw[black,line width=.5pt] (t4)--(t5);
\draw[black,line width=.5pt] (t4)--(t6);
\draw[black,line width=.5pt] (t3)--(t4);
\draw[black,line width=.5pt] (t3)--(t7);
\draw[black,line width=.5pt] (t1)--(t2);
\draw[black,line width=.5pt] (t1)--(t3);
\draw[black,line width=.5pt] (t0)--(t1);
\draw[black,line width=.5pt] (t0)--(t8);
\foreach \v in {0,...,8}{\fill (t\v) circle[radius=.022];}
\end{scope}
\begin{scope}[shift={(a6)}]
\draw[draw=black,fill=white,line width=.55pt] (0,0) circle[radius=.74];
\coordinate (t0) at (0.10281,0.47000);
\coordinate (t1) at (-0.26437,0.23500);
\coordinate (t2) at (-0.47000,0.00000);
\coordinate (t3) at (-0.05875,0.00000);
\coordinate (t4) at (-0.23500,-0.23500);
\coordinate (t5) at (0.11750,-0.23500);
\coordinate (t6) at (0.00000,-0.47000);
\coordinate (t7) at (0.23500,-0.47000);
\coordinate (t8) at (0.47000,0.23500);
\draw[black,line width=.5pt] (t5)--(t6);
\draw[black,line width=.5pt] (t5)--(t7);
\draw[black,line width=.5pt] (t3)--(t4);
\draw[black,line width=.5pt] (t3)--(t5);
\draw[black,line width=.5pt] (t1)--(t2);
\draw[black,line width=.5pt] (t1)--(t3);
\draw[black,line width=.5pt] (t0)--(t1);
\draw[black,line width=.5pt] (t0)--(t8);
\foreach \v in {0,...,8}{\fill (t\v) circle[radius=.022];}
\end{scope}
\begin{scope}[shift={(a7)}]
\draw[draw=black,fill=white,line width=.55pt] (0,0) circle[radius=.74];
\coordinate (t0) at (-0.13219,0.47000);
\coordinate (t1) at (-0.47000,0.23500);
\coordinate (t2) at (0.20563,0.23500);
\coordinate (t3) at (-0.05875,0.00000);
\coordinate (t4) at (-0.23500,-0.23500);
\coordinate (t5) at (0.11750,-0.23500);
\coordinate (t6) at (0.00000,-0.47000);
\coordinate (t7) at (0.23500,-0.47000);
\coordinate (t8) at (0.47000,0.00000);
\draw[black,line width=.5pt] (t5)--(t6);
\draw[black,line width=.5pt] (t5)--(t7);
\draw[black,line width=.5pt] (t3)--(t4);
\draw[black,line width=.5pt] (t3)--(t5);
\draw[black,line width=.5pt] (t2)--(t3);
\draw[black,line width=.5pt] (t2)--(t8);
\draw[black,line width=.5pt] (t0)--(t1);
\draw[black,line width=.5pt] (t0)--(t2);
\foreach \v in {0,...,8}{\fill (t\v) circle[radius=.022];}
\end{scope}
\begin{scope}[shift={(a8)}]
\draw[draw=black,fill=white,line width=.55pt] (0,0) circle[radius=.74];
\coordinate (t0) at (-0.22031,0.47000);
\coordinate (t1) at (-0.47000,0.23500);
\coordinate (t2) at (0.02937,0.23500);
\coordinate (t3) at (-0.23500,0.00000);
\coordinate (t4) at (0.29375,0.00000);
\coordinate (t5) at (0.11750,-0.23500);
\coordinate (t6) at (0.00000,-0.47000);
\coordinate (t7) at (0.23500,-0.47000);
\coordinate (t8) at (0.47000,-0.23500);
\draw[black,line width=.5pt] (t5)--(t6);
\draw[black,line width=.5pt] (t5)--(t7);
\draw[black,line width=.5pt] (t4)--(t5);
\draw[black,line width=.5pt] (t4)--(t8);
\draw[black,line width=.5pt] (t2)--(t3);
\draw[black,line width=.5pt] (t2)--(t4);
\draw[black,line width=.5pt] (t0)--(t1);
\draw[black,line width=.5pt] (t0)--(t2);
\foreach \v in {0,...,8}{\fill (t\v) circle[radius=.022];}
\end{scope}
\end{tikzpicture}
\caption{A cycle of length $9$ in $\mathcal A_6$. Each arrow is a rotation at an associative
triple for the operation $a\cdot b=-a-b+w$ on $\mathbb Z_3$.}
\label{fig:nine-walk}
\end{figure}

Each step of this cycle is a rotation $((U,V),W)\leftrightarrow(U,(V,W))$.
The pairs $(U,W)$ of the nine steps are, in order, $(t,t(tt))$,
$(tt,tt)$, $((tt)t,t)$,
and then $(t,t)$ six times. In every step,
$\lab(U)=\lab(W)$. Hence every step is a rotation at an associative
triple. Since this cycle has odd length, we again obtain a
contradiction.

In every case, we obtain a contradiction. Hence $|S'|\geqslant4$.
Therefore the map $T\mapsto(\lab(T),\beta(T))$ takes values in the set
$S'\times\F$, which has at least $8$ elements.
\end{proof}

Note that
Theorem~\ref{thm:lower} applies only to colourings of the form in its
statement. In particular, Theorem~\ref{thm:lower} does not improve the
lower bound $\chi(\mathcal A_{10})\geqslant4$.

\medskip\noindent\textbf{Acknowledgement.} Research supported by the
Australian Research Council and by the Natural Sciences and Engineering
Research Council of Canada.

\medskip\noindent\textbf{Postscript.} After submitting this paper, the author was informed that Justin Dallant and Jonathan Noel had independently proved that associahedra are 8-colourable.

\end{document}